\documentclass[11pt,reqno]{amsart}
\usepackage{esint}
\usepackage{bbm}
\usepackage{bm}
\usepackage{amssymb}
\usepackage{mathrsfs}
\usepackage{amsfonts}
\usepackage{amsfonts,amssymb,amsmath,amsthm}
\usepackage{url}
\usepackage{enumerate}
\usepackage[pdftex,bookmarksnumbered]{hyperref}
\usepackage{graphicx,xcolor}
\usepackage{pgfplots}
\usepackage{float}
\usepackage{caption}
\usepackage{ulem}
\usepackage{amsmath}
\usepackage{mathtools}

\newtheorem{theorem}{Theorem}[section]

\newtheorem{lem}[theorem]{Lemma}

\newtheorem{prop}[theorem]{Proposition}

\theoremstyle{definition}

\newtheorem{remark}[theorem]{Remark}

\numberwithin{equation}{section}

\newtheorem{ex}{Example}

\newenvironment{proofof}[1]{%
\renewcommand{\proofname}{{\rmfamily\bfseries Proof of Theorem #1}}
\begin{proof}
}{%
\end{proof}
}

\usepackage[left=2.5cm, right=2.5cm, top=2.5cm, bottom=2.5cm]{geometry}

\usepackage{graphicx}%
\usepackage{color}
\usepackage{cite}
\usepackage{hyperref}
\hypersetup{
	colorlinks = true,%
	citecolor = blue,
	filecolor=red,%
	linkcolor = [rgb]{0.65,0.0,0.0},%
	anchorcolor = red,
	pagecolor = red,
	urlcolor= [rgb]{0.65,0.0,0.0}
	linktocpage=true,
	pdfpagelabels=true,
	bookmarksnumbered=true,
}

\makeatletter
\renewcommand{\abstractname}{Abstract} 
\renewenvironment{abstract}{%
	\begin{center}\bfseries \abstractname\end{center} 
	\small 
	\quotation 
}{%
	\endquotation
}
\makeatother

\title{Projectively Flat Spherically Symmetric Sprays of Constant Curvature}
\author{Yajing Gui \and Benling Li}

\thanks{B. Li is supported by Natural Science Foundation of Ningbo (No.~2024J017) and  National Natural Science Foundation of China (No.~12471045). }
\begin{document}
\maketitle
\begin{center}
\vskip -3mm
\small
School of Mathematics and Statistics, Ningbo University, 315211 Ningbo, China\\
{\it libenling@nbu.edu.cn}
\end{center}
\vskip 3mm
\begin{abstract}
In this paper, we study projectively flat spherically symmetric sprays of constant curvature.
We derive a system of equations that is equivalent to the
constant-curvature condition for such sprays and solve it to obtain a complete classification of the projective factor \(p(r,s)\) into four explicit families. We further show that the sprays induced by projectively flat spherically symmetric Finsler metrics of constant flag curvature \(K=0\), \(K=-1\), and \(K=1\) are all special cases of our classification, thereby confirming that the sprays induced by these constant-flag-curvature Finsler metrics are fully contained in our framework.
\end{abstract}

\vskip 5mm
\begin{quotation}

\end{quotation}

\vskip 8mm
\section{Introduction}
Rotational symmetry, equivalently spherical symmetry, i.e., invariance
under the orthogonal group \(O(n)\),
 is one of the most basic structures in geometry. We write
\(|\cdot|\) and \(\langle\cdot,\cdot\rangle\) for the Euclidean norm and
inner product. A fundamental family of examples is provided by the
Riemannian metrics of constant sectional curvature \(\mu\) on a domain
\(\Omega\subseteq\mathbb{R}^n\),
\[
\alpha_\mu
=\frac{\sqrt{|y|^2+\mu\bigl(|x|^2|y|^2-\langle x,y\rangle^2\bigr)}}
{1+\mu|x|^2},
\]
which are rotationally invariant, and whose geodesic sprays are given by
\[
G^i=-\frac{\mu\langle x,y\rangle}{1+\mu|x|^2}\,y^i.
\]
Another classical example is the Schwarzschild solution in general relativity.

Within Finsler geometry, a metric \(F\) on
\(\Omega\subseteq\mathbb{R}^n\) is said to be \textit{spherically
symmetric} if \(F(Ux,Uy)=F(x,y)\) for every \(U\in O(n)\) and every
\(y\in T_x\Omega\). The systematic study of such metrics goes back to
McCarthy--Rutz~\cite{McCarthyRutz1993}, who derived the general form of
four-dimensional spherically symmetric Finsler spaces. In arbitrary
dimension, Zhou~\cite{Zhou2012b} established that every spherically
symmetric Finsler metric is of the form
\(F=\varphi(|x|,|y|,\langle x,y\rangle)\); he also classified the
projectively flat ones of constant flag curvature~\cite{Zhou2012a}.
The scalar-curvature case was treated by Huang--Mo~\cite{HuangMo2012},
who obtained the governing equations and produced non-projectively flat
examples. There is a substantial body of related work; see,
e.g.,~\cite{AsmaZoltan2024,Elgendi2024,IlmavirtaMonkkonen2023,
KhaniMoghaddam2025,Li2020,LiMoYu2019,SolorzanoLeandro2021,
TayebiNajafi2025}. These results indicate that rotational symmetry is a
natural and fruitful notion beyond the metric setting,
which motivates considering it for sprays, defined below.

Let \(M\) be a smooth \(n\)-dimensional manifold with \(n\ge2\). A
\textit{spray} on \(M\) is a smooth vector field \(G\) on
\(TM\setminus\{0\}\) which, in local coordinates \((x^i,y^i)\), takes
the form
\[
G=y^i\frac{\partial}{\partial x^i}-2G^i\frac{\partial}{\partial y^i},
\]
where the coefficients \(G^i=G^i(x,y)\) are smooth on
\(TM\setminus\{0\}\) and positively homogeneous of degree two in \(y\),
i.e.,
\[
G^i(x,\lambda y)=\lambda^2G^i(x,y),\qquad \forall\,\lambda>0.
\]
Every Finsler metric gives rise to a spray, but the converse fails:
there are many sprays that are not induced by any Finsler metric.

Since rotational invariance is naturally formulated on domains in
\(\mathbb{R}^n\), we consider sprays on a domain
\(\Omega\subseteq\mathbb{R}^n\). Such a spray is said to be
\textit{spherically symmetric} if it is invariant under all rotations,
that is,
\[
G(Ux,Uy)=G(x,y),\qquad \forall\,U\in O(n).
\]
Recently, Gui--Li~\cite{GuiLi2026} showed that a spray on \(\Omega\) is
spherically symmetric if and only if its coefficients admit the
representation
\begin{equation}\label{eq:projective_spherical_form}
G^i(x,y)=|y|\,\alpha(r,s)\,y^i+|y|^2\,\beta(r,s)\,x^i
\quad(y\neq0),\qquad G^i(x,0)=0,
\end{equation}
where \(r=|x|^2\), \(s=\langle x,y\rangle/|y|\), and \(\alpha,\beta\)
are smooth functions depending only on \(r\) and \(s\).

Projective flatness is a classical topic in Finsler geometry, closely
related to classification problems and to Hilbert's fourth problem.
Recall that a spray on a domain \(\Omega\subseteq\mathbb{R}^n\) is
called \textit{projectively flat} if all its geodesics are straight
lines, i.e.,
\[
G^i=P(x,y)y^i,
\]
where \(P(x,y)\), positively homogeneous of degree one in \(y\), is
called the \textit{projective factor} of \(G\). Studying projective
flatness in the spray setting is therefore both natural and of
independent interest. In connection with Hilbert's fourth problem,
Buc\u{a}taru--Muzsnay~\cite{BucataruMuzsnay2014} characterized when
projectively flat sprays are metrizable by Finsler metrics.
Yang~\cite{Yang2011} studied a class of projectively flat sprays and
showed that there exist many isotropic sprays that cannot be induced by
any (even singular) Finsler metric. For a spray that is both projectively flat
and spherically symmetric, comparing $G^i=P(x,y)y^i$ with the Gui--Li
representation above shows that the coefficient of $x^i$ must vanish,
so that
\[
G^i=|y|\,p(r,s)\,y^i,\qquad
P(x,y)=|y|\,p(r,s),
\]
where $r$ and $s$ are as above.

We now turn to curvature. In 1926, Berwald~\cite{Berwald1927} extended
the notion of the Riemann curvature tensor to sprays. A spray $G$ is
said to be of \textit{scalar curvature} if its Riemann curvature tensor
$R^i_{\ k}$ takes the form
\[
R^i_{\ k}=R\delta^i_k-\tau_ky^i,
\]
where $R=R(x,y)$, $\tau_k=\tau_k(x,y)$, and $R=\tau_ky^k$. If in
addition $\tau_k=\frac{1}{2}R_{y^k}$, then $G$ is said to be of
\textit{isotropic curvature}. This notion was introduced by
Li--Shen~\cite{LiShen2018}, who showed that a spray of scalar curvature
is isotropic exactly when the non-Riemannian $\chi$-curvature vanishes;
they also determined the structure of projectively flat sprays with
nonzero isotropic curvature. A further generalization, weakly isotropic
curvature, was introduced by Lou--Li~\cite{LouLi2022}.

In 2023, Yang~\cite{Yang2023} introduced sprays of \textit{constant
curvature}, defined by
\[
\tau_{i;k}=0,
\]
where $\tau_{i;k}$ denotes the horizontal covariant derivative of the
component $\tau_i$ with respect to $x^k$. He pointed out that a spray of
isotropic curvature need not be of constant curvature, even in dimension
$n\ge3$. At the level of Finsler metrics, Li~\cite{Li2014} obtained
complete classifications of locally projectively flat Finsler metrics
with constant flag curvature $K=0,-1,1$. Further results on spherically
symmetric Finsler metrics of constant curvature can be found
in~\cite{BucataruMuzsnay2013,CaiQiuWang2020,SevimShenUlgen2015,XuLi2013}.

In this paper, we study projectively flat spherically symmetric sprays
of constant curvature. By \eqref{eq:projective_spherical_form}, such a
spray has coefficients \(G^i=|y|p(r,s)y^i\). We give a complete
classification of the projective factor \(p(r,s)\). Our main result is
the following.
\begin{theorem}\label{maintheorem1}
Let \(G\) be a projectively flat spherically symmetric spray on a
domain \(\Omega\subseteq\mathbb{R}^n\) with \(n\ge2\) and
 \(0\in\Omega\), with
\(G^i=|y|p(r,s)y^i\) as in \eqref{eq:projective_spherical_form}.
Then \(G\) has constant curvature if and
only if \(p(r,s)\) is given by one of the following forms, each
understood on its natural domain:
\begin{enumerate}[\rm (i)]
\item
\begin{equation}\label{sol1}
p(r,s)\equiv0;
\end{equation}
\item
\begin{equation}\label{sol2}
p(r,s)=\frac{s+\varepsilon\sqrt{c_1-r+s^2}}{2(c_1-r)},
\end{equation}
where \(\varepsilon=\pm1\) and \(c_1\) is a positive constant;
\item
\begin{equation}\label{sol3}
p(r,s)=-\frac{s}{r+c_2},
\end{equation}
where \(c_2\) is a nonzero constant;
\item
\begin{equation}\label{sol4}
p(r,s)=\frac{\eta(\delta s^{2}-Q)\sqrt{\delta Q-r+s^{2}-c_{3}}}
{\sqrt2(r^2+2c_{3}r+c_{4})}
-\frac{(r+c_{3})s}{r^2+2c_{3}r+c_{4}},
\end{equation}
where \(\eta=\pm1\), \(\delta=\pm1\), and
\(Q=\sqrt{(r-s^2)^2+2c_3(r-s^2)+c_4}\), with \(c_3\) and \(c_4\)
constants such that \eqref{sol4} is well-defined; that is,
\[
\delta Q-r+s^{2}-c_{3}\ge0
\quad \text{and}\quad
(r-s^2)^2+2c_3(r-s^2)+c_4>0.
\]
Concrete admissible choices of \(c_3,c_4\) are given in
Remark~\ref{remarkthm}.
\end{enumerate}
\end{theorem}

Here the natural domain means the set of \((r,s)\) for which all
radicals are real and all denominators are nonzero, unless a more
specific admissible range is stated in (iv).

Theorem~\ref{maintheorem1} shows that the constant-curvature condition
singles out exactly four families of projective factors among
projectively flat spherically symmetric sprays. The first three
families are elementary and explicit, while the fourth depends on two
parameters $c_3,c_4$ and captures the genuinely non-Riemannian
phenomena. In Section~\ref{4}, we compare these four families with the
classification of projectively flat spherically symmetric Finsler
metrics of constant flag curvature $K=0,-1,1$. Such a classification
was first obtained by Zhou~\cite{Zhou2012a}, and a complete explicit
form was later given by Li~\cite{Li2014} (recalled in
Lemma~\ref{classification of K}). The sprays induced by those metrics
correspond to specific parameter choices in (i)--(iv). Since this
classification is complete, every other parameter choice in
Theorem~\ref{maintheorem1} yields a constant-curvature spray that does
not arise from any of the Finsler metrics classified in
Lemma~\ref{classification of K}. In this sense, within the projectively
flat spherically symmetric setting, the spray classification of
Theorem~\ref{maintheorem1} is strictly larger than its Finslerian
counterpart, and it exhibits explicitly the constant-curvature sprays
strictly larger than the class of projectively flat spherically symmetric Finsler metrics of constant flag curvature.

\section{Preliminaries}\label{2}
In this section, we collect the definitions and known results used in
the proof of Theorem~\ref{maintheorem1}. Standard references for spray
geometry and Finsler geometry are the monograph of Shen~\cite{ShenBook}
and the survey of Li--Shen~\cite{LiShen2024Survey}. Throughout this paper, we use the Einstein summation convention:
repeated upper and lower indices are summed from $1$ to $n$.

Let $M$ be an $n$-dimensional smooth manifold with $n\geq2$. A \textit{spray} $G$ on  $M$ is a smooth vector field on $TM\setminus\{0\}$, which in local coordinates \((x^i,y^i)\) takes the form
\begin{equation*}
	G=y^i\frac{\partial}{\partial x^i}-2G^i\frac{\partial}{\partial y^i},
\end{equation*}
where the coefficients \(G^i = G^i(x, y)\)  are smooth on \(TM \setminus \{0\}\) and positively homogeneous of degree two in $y$. The \textit{geodesic coefficients} of $G$ can be expressed by a system of second order differential equations
\begin{equation*}
	\frac{d^2x^i}{dt^2}+2G^i(x,\frac{dx}{dt})=0.
\end{equation*}

A spray on a domain \(\Omega\subseteq\mathbb{R}^n\) is called \textit{projectively flat} if all its geodesics are straight lines, i.e.
\[
G^i = P(x,y)y^i,
\]
where \(P(x,y)\) is homogeneous of degree one in \(y\) and  called {\it projective factor} of $G$.

A spray on a domain $\Omega$ is called \textit{spherically symmetric} if it is invariant under all rotations:
\begin{equation*}
	G(Ux,Uy)=G(x,y),\qquad \forall\, U\in O(n).
\end{equation*}
Under the spherically symmetric condition, the geodesic coefficients of $G$ depend only on \(|x|,|y|,\langle x,y\rangle\). Let \(r=|x|^2\), \(s=\langle x,y\rangle/|y|\). For a projectively flat spherically symmetric spray, its projective factor has the form
\[
P(x,y)=|y|\,p(r,s).
\]

The \textit{Riemann curvature tensor} $R^i_{\ k}$ is a linear transformation on the tangent space, which in local coordinates can be expressed as
\begin{equation*}
	R^i_{\ k}=2  G^i_{x^k}- G^i_{x^j y^k} y^j+2 G^j G^i_{y^j y^k}- G^i_{y^j} G^j_{y^k}.
\end{equation*}
The following lemma is about the Riemann curvature tensor of a projectively flat spherically symmetric spray.
\begin{lem}\cite{GuiLi2026}\label{R}
	Let $G$ be a projectively flat spherically symmetric spray with $G^i = |y| p(r,s) y^i$, where $r=|x|^2$ and $s=\langle x,y\rangle/|y|$, then the Riemann curvature tensor of such a spray is given by
	\begin{align*}
		\begin{split}
			R^i_{\ k}&=|y|^2(p^2-2sp_r-p_s)\delta^i_k\\
			&\quad -\Big\{|y|(pp_s+2sp_{rs}-4p_r+p_{ss})x_k
			-\big[(sp+1)p_s+2s(sp_{rs}-p_r)+sp_{ss}-p^2\big]y_k\Big\}y^i,
		\end{split}
	\end{align*}
	and
	\begin{equation*}
		R=|y|^2(p^2-2sp_r-p_s).
	\end{equation*}
\end{lem}
A spray $G$ is said to be of \textit{scalar curvature} if its Riemann curvature tensor $R^i_{\ k}$ satisfies
\[
R^i_{\ k}=R\delta^i_k-\tau_ky^i,
\]
where $R=R(x,y)$, $\tau_k=\tau_k(x,y)$ with $R=\tau_ky^k$. If $\tau_k=\frac12R_{y^k}$,
the spray $G$ is of \textit{isotropic curvature}.
The following lemma is about the relation that \(p(r,s)\) must satisfy for a projectively flat spherically symmetric spray with isotropic curvature.

\begin{lem}\cite{GuiLi2026}\label{isotropic}
	Let \(G\) be a projectively flat spherically symmetric spray with \(G^i=|y|p(r,s)y^i\), where \(r=|x|^2\), \(s=\langle x,y\rangle/|y|\). Then \(G\) has  isotropic curvature if and only if
	\begin{equation}\label{isotropic eq}
		p_{ss}=2p_r-2sp_{rs}.
	\end{equation}
\end{lem}

As a special case of isotropic curvature, a spray $G$ is called of {\it constant curvature} if
\begin{equation}\label{taui_k}
\tau_{i;k}=0.
\end{equation}
 Here, \(\tau_{i;k}\) denotes the horizontal covariant derivative of \(\tau_{i}\) with respect to \(x^k\). The explicit operation is $\frac{\delta}{\delta x^i}=\frac{\partial}{\partial x^i}-N^j_i\frac{\partial}{\partial y^j}$, where $N^j_i=\frac{\partial G^j}{\partial y^i}$.

The following proposition concerns sprays of scalar curvature, which will help us solve the equations for projectively flat spherically symmetric sprays of constant curvature later.
\begin{prop}\cite{Yang2023}\label{Prop2.3}
Let $G$ be of scalar curvature, then
\[
R_{;i}=0\iff\tau_{i;k} =0.
\]
\end{prop}
Since a projectively flat spray is of scalar curvature (cf.~\cite{ShenBook}), Proposition~\ref{Prop2.3} applies.
By the above proposition and \eqref{taui_k}, we see that a spray $G$ is of constant curvature if and only if
\begin{equation}\label{R_xi=0}
    R_{;i}=R_{x^i}-N^j_iR_{y^j}=0.
\end{equation}

A \textit{Finsler metric} on $M$ is a nonnegative function $F$: $TM\longrightarrow [0,  +\infty)$ satisfying the following conditions: (i) $F$ is smooth on $TM_0:=TM \setminus \{0\}$; (ii) $F(x,\lambda y)=\lambda F(x,y)$, $  \forall \lambda>0$, i.e. $F$ is a positively homogeneous function of degree one in $y$; (iii) $\forall y\neq0$, the matrix $g_{ij}(x,y):=\frac{1}{2}[F^2]_{y^iy^j}$ is positive definite. Every Finsler metric on $M$ induces a spray by
\begin{equation*}
	G^i=\frac{1}{4}g^{ij}\{[F^2]_{y^jx^k}y^k-[F^2]_{x^j}\},
\end{equation*}
where $g^{ij}$ is the inverse matrix of $g_{ij}$.

Let $P=\operatorname{span}\{y,u\}$ be a 2-dimensional plane in $T_xM$, where $y,u \in T_xM \setminus \{0\}$. The \textit{flag curvature} $K(P,y)$ with $y$ as the flagpole is defined as
\begin{equation*}
	K(P,y) = \frac{g_y(u, R_y(u))}{g_y(y,y)g_y(u,u) - \left[g_y(y,u)\right]^2},
\end{equation*}
where $g_y = g_{ij}(x,y)dx^i \otimes dx^j$ and $R_y(u)= R^i_{\ k}dx^k \otimes \frac{\partial}{\partial x^i}$.
A Finsler metric has \textit{scalar flag curvature} if $K(P,y)$
 depends only on the base point \(x\) and the flagpole \(y\), and \textit{constant flag curvature} if $K\equiv\lambda$ is a constant.

If all geodesics of \(F\) are straight lines, then \(F\) is said to be \textit{projectively flat}; its geodesic coefficients $G^i$ satisfy
\begin{equation}\label{P_F}
G^i = P(x, y)y^i,\qquad
P = \frac{F_{x^k}y^k}{2F}
\end{equation}
where $P$ is the \textit{projective factor} of \(F\).
 The following lemma gives a classification of projectively flat spherically symmetric Finsler metrics with constant flag curvature.
\begin{lem}\cite{Li2014}\label{classification of K}
Let $F$ be a projectively flat spherically symmetric Finsler metric on a convex domain \(\Omega\subseteq \mathbb{R}^n\).
If $F$ is of constant flag curvature $K \in \{0, -1, 1\}$, the classification is as follows:
\begin{enumerate} [\rm (i)]
\item $K=0$ if and only if
\begin{equation}\label{K=0_F}
F=|y|\quad \text{or} \quad
F = \frac{|y|^4}{z\bigl(c\langle x,y\rangle \pm z\bigr)^2},
\end{equation}
where $z = \sqrt{(1 - c^2|x|^2)|y|^2 + c^2\langle x,y\rangle^2}$, $c$ is a nonzero constant;
\item $K=-1$ if and only if
\begin{align}
F&=\frac12\left\{
\frac{(c+1)^2\langle x,y\rangle+\operatorname{sgn}(c+1)\sqrt{(c+1)^2\bigl(1-(c+1)^2|x|^2\bigr)|y|^2+(c+1)^4\langle x,y\rangle^2}}{1-(c+1)^2|x|^2}\right.\nonumber\\
&\quad\left.-\frac{(c-1)^2\langle x,y\rangle+\operatorname{sgn}(c-1)\sqrt{(c-1)^2\bigl(1-(c-1)^2|x|^2\bigr)|y|^2+(c-1)^4\langle x,y\rangle^2}}{1-(c-1)^2|x|^2}
\right\},\label{K=-1_F}
\end{align}
where $c$ is a constant.
\item
$K=1$ if and only if
\begin{equation}\label{K=1_F}
F=\operatorname{Im}\left[
\frac{(c+i)^2\langle x,y\rangle+\sqrt{(c+i)^2\bigl(1-(c+i)^2|x|^2\bigr)|y|^2+(c+i)^4\langle x,y\rangle^2}}
{1-(c+i)^2|x|^2}
\right],
\end{equation}
where $c$ is a constant. Obviously, when $c=0$
\[F_{c=0} = \frac{\sqrt{(1+|x|^2)|y|^2 - \langle x,y\rangle^2}}{1+|x|^2}.\]

\end{enumerate}
\end{lem}

\vskip 10mm
\section{Projectively flat spherically symmetric sprays of constant curvature}\label{3}
In this section, we present an equivalent system of equations for a
projectively flat spherically symmetric spray of constant curvature,
and provide explicit solutions. The proof of the main theorem (Theorem~\ref{maintheorem1}) is split into two lemmas and a final classification step.

Our first lemma gives an equivalent PDE system for a projectively flat spherically symmetric spray $G$ of constant curvature.
\begin{lem}\label{lem:constsys}
	Let \(G\) be a projectively flat spherically symmetric spray with \(G^i=|y|p(r,s)y^i\), where \(r=|x|^2\), \(s=\langle x,y\rangle/|y|\). Then \(G\) has constant curvature if and only if \(p(r,s)\) satisfies the following system:
	\begin{align}\label{const}
		\begin{cases}
			-4sp_{rr}+pp_{ss}+2(sp-1)p_{rs}+2(2sp_{s}+3p)p_r-2(2p^2-p_s)p_s=0,\\[3pt]
			-(sp+1)(p_{ss}+2sp_{rs})-2(2s^{2}p_{s}-3sp+1)p_r+2(2sp^2-sp_s+3p)p_s-4p^{3}=0.
		\end{cases}
	\end{align}
Moreover, if \(p\) satisfies \eqref{const}, then it necessarily satisfies the isotropic curvature condition \eqref{isotropic eq}.
Substituting \eqref{isotropic eq} into \eqref{const} to eliminate
\(p_{ss}\), and retaining \eqref{isotropic eq} as an additional
equation, yields the equivalent system
	\begin{align}\label{const1}
		\begin{cases}
			-4p^2p_s+4s p_r p_s-4s p_{rr}+8p p_r+2p_s^2-2p_{rs}=0,\\
			-2s p_s^2+(-4s^2p_r+4sp^2+6p)p_s+4s p p_r-4p^3-4p_r=0,\\
			p_{ss}-2p_r+2sp_{rs}=0.
		\end{cases}
	\end{align}
\end{lem}
\begin{proof}
By \eqref{R_xi=0}, a projectively flat spherically symmetric spray $G$ is of constant curvature if and only if $R_{;i}=0$. In this case, by Lemma \ref{R}, we have \[ R=|y|^2(p^2-2sp_r-p_s).\]
 To compute $R_{;i}$, we need the derivatives of $p$ with respect to $x^i$ and $y^j$. Since $p=p(r,s)$ with $r=|x|^2, s=\langle x,y\rangle/|y|$, by the chain rule, we have
	\begin{equation*}
		p_{x^i}=p_r r_{x^i}+p_s s_{x^i},\qquad p_{y^j}=p_s s_{y^j},
	\end{equation*}
	where
	\begin{equation*}
		r_{x^i}=2x_i,\quad s_{x^i}=\frac{y_i}{|y|},\quad s_{y^j}=\frac{|y|^2x_j-\langle x,y\rangle y_j}{|y|^3},
	\end{equation*}
and analogous formulas hold for the derivatives of $p_r$ and $p_s$. Applying these to compute $N^j_i=G^j_{y^i}$, $R_{x^i}$, and $R_{y^j}$, we obtain after simplification
\begin{align*}
N^j_i&=\frac{y^j}{|y|}\left[|y| p_sx_i+(p-sp_s)y_i\right]+|y| p\delta^j_i,\\
R_{x^i}&=|y|^2(4pp_r-4sp_{rr}-2p_{rs})x_i+|y|(2pp_s-2sp_{rs}-2p_r-p_{ss})y_i,\\
R_{y^j}&=|y|(2pp_s-2sp_{rs}-2p_r-p_{ss})x_j+\left[s(2sp_{rs}-2pp_s+p_{ss}-2p_r)+2(p^2-p_s)\right]y_j.
\end{align*}
Substituting the above equations back into $R_{;i}$ in \eqref{R_xi=0}, we get
\begin{align*}
R_{;i} & =|y|^2\left[-4sp_{rr}+pp_{ss}+2(sp-1)p_{rs}+2(2sp_{s}+3p)p_r-2(2p^2-p_s)p_s\right]x_i\\
&\quad +|y|\left[-(sp+1)(p_{ss}+2sp_{rs})-2(2s^{2}p_{s}-3sp+1)p_r+2(2sp^2-sp_s+3p)p_s-4p^{3}\right]y_i.
\end{align*}
Since \(R_{;i}=0\) for all \(y\ne0\), taking \(y\) not parallel to \(x\)
makes \(x_i\) and \(y_i\) linearly independent; hence both coefficients
vanish. The cases where \(y\) is parallel to \(x\), and where \(x=0\), follow
by continuity. Therefore
system \eqref{const} holds. This proves the necessity.

Conversely, if \(p\) satisfies \eqref{const}, then the coefficients of
\(x_i\) and \(y_i\) in the above expression for \(R_{;i}\) vanish
identically, and hence \(R_{;i}=0\). By \eqref{R_xi=0}, \(G\) is of
constant curvature. This proves the sufficiency.

It remains to justify the equivalence between \eqref{const} and
\eqref{const1}.
Since a spray of constant curvature is necessarily of isotropic
curvature, and since \(G\) has constant curvature exactly when \(p\)
satisfies \eqref{const}, it follows that \eqref{const} implies
\eqref{isotropic eq}. Thus \eqref{const} is equivalent to the system
formed by \eqref{const} and \eqref{isotropic eq}. Substituting the
expression for \(p_{ss}\) from \eqref{isotropic eq} into \eqref{const},
and retaining \eqref{isotropic eq} as an
additional equation, yields \eqref{const1}.
\end{proof}

To solve the equivalent system of equations \eqref{const1}, motivated by the solution of the isotropic curvature condition in \cite{GuiLi2026}, we perform the change of variables \(u=r-s^2, v=s\); in this case, $p(r,s)=\tilde{p}(u,v)$. By the chain rule, \eqref{const1} is transformed into
\begin{align}
	&-\tilde{p}_{uv}+(4v\tilde{p}^2-2v\tilde{p}_v+4\tilde{p})\tilde{p}_u-(2\tilde{p}^2-\tilde{p}_v)\tilde{p}_v=0,\label{treq1}\\
	&(2v^2\tilde{p}_v-4v^2\tilde{p}^2-4v\tilde{p}-2)\tilde{p}_u+(2v\tilde{p}^2+3\tilde{p}-v\tilde{p}_v)\tilde{p}_v-2\tilde{p}^3=0,\label{treq2}\\
	&-2v\tilde{p}_{uv}-4\tilde{p}_u+\tilde{p}_{vv}=0.\label{treq3}
\end{align}
Since the change of variables \((r,s)\mapsto(u,v)=(r-s^2,s)\) has
inverse \(r=u+v^2\), \(s=v\), and its Jacobian matrix is invertible,
system \eqref{const1} is equivalent to \eqref{treq1}--\eqref{treq3}.
Before solving \eqref{treq1}--\eqref{treq3}, we prove the following
preliminary lemma.
\begin{lem}\label{ABform}
Let \(G\) be a projectively flat spherically symmetric spray with \(G^i=|y|p(r,s)y^i\), where \(r=|x|^2\), \(s=\langle x,y\rangle/|y|\).
 Set \(u=r-s^2\), \(v=s\), and
\(\tilde p(u,v)=p(r,s)\). If \(\tilde p\) satisfies
\eqref{treq1}--\eqref{treq3}, then
	\begin{equation}\label{prepform}
	\tilde{p}(u,v)=-\frac{2 A v+B }{2A v^2+2B v+2},
	\end{equation}
where \(A=A(u)\) and \(B=B(u)\) satisfy
	\begin{align}
		-B^2 A+2A^2+2A'&=0, \label{ABeq2}\\
		-B^3+3AB+2B'&=0. \label{ABeq3}
	\end{align}
\end{lem}

\begin{proof}
We first derive two consequences of \eqref{treq1}--\eqref{treq3}.
Subtracting \(2v\times\eqref{treq1}\) from \eqref{treq3} eliminates
\(\tilde{p}_{uv}\) and yields
\begin{equation}\label{eq31}
	(4v^2\tilde{p}_v-8v^2\tilde{p}^2-8v\tilde{p}-4)\tilde{p}_u
	+(4v\tilde{p}^2-2v\tilde{p}_v)\tilde{p}_v+\tilde{p}_{vv}=0.
\end{equation}
The combination \(2\times\eqref{treq2}-\eqref{eq31}\) yields
\begin{equation}\label{u_betabeta}
	6\tilde{p}\tilde{p}_v-4\tilde{p}^3-\tilde{p}_{vv}=0.
\end{equation}

Next we reduce \eqref{u_betabeta} to a linear equation. Fix an
admissible \(u>0\) and a point \(v_0\) in the interior of the
corresponding \(v\)-interval. Define
\begin{equation}\label{pD}
	h(u,v):=\exp\left(-2\int_{v_0}^{v}\tilde p(u,\tau)\,d\tau\right),
\end{equation}
where the integral is taken in the \(v\)-direction with \(u\) fixed
as a parameter. Since \(\tilde p\) is smooth on the domain under
consideration, \(h\) is well defined and smooth on the same domain,
and \(h(u,v_0)=1\); in particular, \(h(u,v)\neq0\). Moreover,
differentiating \eqref{pD} with respect to \(v\) gives
\[
h_v=-2\tilde p\,h,
\]
equivalently
\begin{equation}\label{pD2}
	\tilde p(u,v)=-\frac{h_v(u,v)}{2h(u,v)}.
\end{equation}

A direct computation from \eqref{pD2} gives the identity
\begin{equation}\label{identity}
	6\tilde p\tilde p_v-4\tilde p^3-\tilde p_{vv}
	=\frac{h_{vvv}}{2h}.
\end{equation}
Indeed, setting \(w:=h_v/h\), so that \(\tilde p=-w/2\), one has
\[
6\tilde p\tilde p_v-4\tilde p^3-\tilde p_{vv}
=\frac12\left(3ww_v+w^3+w_{vv}\right),
\]
and substituting
\[
w_v=\frac{hh_{vv}-h_v^2}{h^2},\qquad
w_{vv}=\frac{h^2h_{vvv}-3hh_vh_{vv}+2h_v^3}{h^3}
\]
yields \eqref{identity}. Since \(h\neq0\), equations
\eqref{u_betabeta} and \eqref{identity} together are equivalent to
\[
h_{vvv}=0.
\]
Hence, for functions \(f_2,f_1,f_0\) of \(u\) alone,
\begin{equation*}
	h(u,v)=f_2(u)v^2+f_1(u)v+f_0(u).
\end{equation*}

We claim that \(f_0(u)\neq0\) for all admissible \(u>0\).
Indeed, for each such \(u\), the domain contains an interval around
\(v=0\) on which \(h(u,\cdot)\) is nowhere zero; evaluating at
\(v=0\) gives
\[
h(u,0)=f_0(u)\neq0.
\]
The boundary case \(u=0\) follows by continuity, whenever \(u=0\)
belongs to the closure of the domain.

Since \(\tilde p=-h_v/(2h)\) is invariant under
\(h\mapsto\lambda(u)h\) for any nonzero function \(\lambda\) of
\(u\), we may normalize \(h\) by taking \(\lambda(u)=1/f_0(u)\).
Equivalently, we may assume \(f_0(u)\equiv1\) and write
\begin{equation}\label{D2}
	h(u,v)=A v^2+B v+1,
\end{equation}
where
\[
A =\frac{f_2(u)}{f_0(u)},\qquad
B =\frac{f_1(u)}{f_0(u)}.
\]
Combining \eqref{pD2} with \eqref{D2} yields
\begin{equation}\label{pre_u}
	\tilde p(u,v)
	=-\frac{2A v+B }{2A v^2+2B v+2}.
\end{equation}

We now determine \(A\) and \(B\). Substituting \eqref{pre_u} into
\eqref{treq2} yields
\begin{equation}\label{num}
	\frac{(-BA^2+2BA'-2AB')v^2+(-2B^2A+4A^2+4A')v
	-B^3+3AB+2B'}{2(Av^2+Bv+1)^3}=0,
\end{equation}
where \(A'\) and \(B'\) denote the derivatives of \(A\) and \(B\)
with respect to \(u\). Since \(v\) is arbitrary and the numerator
in \eqref{num} is a polynomial in \(v\) vanishing on a nonempty
open set, all of its coefficients must vanish:
\begin{align}
	-BA^2+2BA'-2AB'&=0,\label{AB1}\\
	-B^2A+2A^2+2A'&=0,\label{AB2}\\
	-B^3+3AB+2B'&=0.\label{AB3}
\end{align}
From \eqref{AB2} and \eqref{AB3} we obtain
\begin{equation}\label{A_uB_u}
	A'=\frac12\left(B^2A-2A^2\right),\qquad
	B'=\frac12\left(B^3-3AB\right).
\end{equation}
Substituting \eqref{A_uB_u} into the left-hand side of \eqref{AB1}
gives
\[
-BA^2+2B\cdot\frac{B^2A-2A^2}{2}
-2A\cdot\frac{B^3-3AB}{2}
=-BA^2+B^3A-2A^2B-AB^3+3A^2B=0,
\]
so \eqref{AB1} is a consequence of \eqref{AB2} and \eqref{AB3}.
Therefore system \eqref{AB1}--\eqref{AB3} is equivalent to system
\eqref{AB2}--\eqref{AB3}. This completes the proof.
\end{proof}

\begin{proofof}{\ref{maintheorem1}}
By Lemmas \ref{lem:constsys} and \ref{ABform}, it remains to solve the
ODE system \eqref{ABeq2}--\eqref{ABeq3} and to verify that each
solution, substituted into \eqref{prepform}, yields a function
\(\tilde p\) satisfying \eqref{const1}; by Lemma \ref{lem:constsys},
the corresponding spray then has constant curvature.

Since each of \(A\) and \(B\) is either identically zero or not, the
following four cases are exhaustive and mutually exclusive. In
Cases 2--4 we solve the system on any open interval where the
relevant functions do not vanish; the remaining points are obtained
by continuity. In each case, the verification is carried out by
direct substitution into \eqref{const1}.

We consider four distinct cases:
\textbf{Case 1:} \(A=B\equiv0\);
\textbf{Case 2:} \(A\equiv0,\ B\not\equiv 0\);
\textbf{Case 3:} \(A\not\equiv 0,\ B\equiv0\);
and \textbf{Case 4:} \(A\not\equiv 0,\ B\not\equiv 0\).

\medskip

\textbf{Case 1:}
If \(A=B\equiv0\), then \eqref{prepform} gives
\(p(r,s)\equiv0\). Conversely, if \(p\equiv0\), then \(p\) satisfies
\eqref{const1}; by the sufficiency part of Lemma \ref{lem:constsys},
\(R_{;i}=0\), and consequently \(G\) has constant curvature.

\medskip
\textbf{Case 2:} If $A\equiv0$, \eqref{ABeq2} is automatically satisfied, and \eqref{ABeq3} becomes
\begin{equation*}
	\frac{dB}{du}=\frac{1}{2}B^3,
\end{equation*}
hence
\begin{equation*}
	B=\frac{\varepsilon}{\sqrt{c_1-u}},
\end{equation*}
where $\varepsilon=\pm1$ and $c_1>0$. In this case,
\begin{equation*}
\tilde{p}(u,v)=-\frac{B}{2Bv+2}=-\frac{1}{2(v+\varepsilon\sqrt{c_1-u})}.
\end{equation*}
Replacing \(\varepsilon\) by \(-\varepsilon\) (which is harmless since
\(\varepsilon=\pm1\)) and using \(p(r,s)=\tilde p(u,v)\), we obtain
\begin{equation}\label{case2}
	p(r,s)=\frac{s+\varepsilon\sqrt{c_1-r+s^2}}{2(c_1-r)}.
\end{equation}
Conversely, if \(p(r,s)\) is given by \eqref{case2}, then its partial
derivatives of all orders are
\begin{align*}
p_r&=\frac{2s\sqrt{s^2-r+c_1}+\varepsilon\left(2s^2-r+c_1\right)}
{4\left(r-c_1\right)^2\sqrt{s^2-r+c_1}},\qquad p_s=-\frac{\varepsilon s+\sqrt{s^2-r+c_1}}
{2\left(r-c_1\right)\sqrt{s^2-r+c_1}},\\
p_{ss}&=\frac{\varepsilon}
{2\left(s^2-r+c_1\right)^{3/2}},\qquad p_{rs}=\frac{2\left(s^2-r+c_1\right)^{3/2}
-\varepsilon s\left(3r-2s^2-3c_1\right)}
{4\left(r-c_1\right)^2\left(s^2-r+c_1\right)^{3/2}},\\
p_{rr}&=\frac{8s\left(r-s^2-c_1\right)\sqrt{s^2-r+c_1}
-\varepsilon\left(8s^4-12s^2\left(r-c_1\right)
+3\left(r-c_1\right)^2\right)}
{8\left(r-c_1\right)^3\left(s^2-r+c_1\right)^{3/2}}.
\end{align*}
Substituting these derivatives into \eqref{const1} verifies that
\eqref{const1} is satisfied. Therefore, if the function \(p\) in the
projective factor \(P=|y|p(r,s)\) is of the form \eqref{case2}, then
the spray has constant curvature.

\medskip

\textbf{Case 3:} If $B\equiv0$, \eqref{ABeq3} is automatically satisfied, and \eqref{ABeq2} becomes
\begin{equation*}
	\frac{dA}{du}=-A^2.
\end{equation*}
Solving this ODE gives
\begin{equation*}
	A =\frac{1}{u+c_2},
\end{equation*}
for some constant \(c_2\). Since \(0\in\Omega\), the point \(x=0\)
belongs to \(\Omega\), so the projective factor must be smooth on
\(T_{0}\Omega\setminus\{0\}\). If \(c_2=0\), then
\begin{equation*}
	\tilde p(u,v)=-\frac{v}{v^2+u}=-\frac{s}{r},
\end{equation*}
which has no smooth extension to \(r=0\); indeed, along \(x\to0\)
with \(\langle x,y\rangle\ne0\), one has
\(s/r=\langle x,y\rangle/(|y||x|^2)\to\infty\) in magnitude.
Thus \(c_2=0\) is excluded by \(0\in\Omega\), and we must have
\(c_2\ne0\). In this case,
\begin{equation*}
	\tilde{p}(u,v)=-\frac{A v}{A v^2+1}=-\frac{v}{v^2+u+c_2}.
\end{equation*}
Using $p(r,s)=\tilde{p}(u,v)$, we have
\begin{equation}\label{case3}
	p(r,s)=-\frac{s}{r+c_2}.
\end{equation}
 Conversely, if $p(r,s)$ satisfies \eqref{case3}, then its partial derivatives of all orders are
\begin{align*}
p_r=\frac{s}{(r+c_2)^2},\quad p_s=-\frac{1}{r+c_2},\quad
p_{ss}=0,\quad
p_{rs}=\frac{1}{(r+c_2)^2},\quad p_{rr}=-\frac{2s}{(r+c_2)^3}.
\end{align*}
Substituting these into \eqref{const1}, we see that all three
equations in \eqref{const1} are satisfied.
Therefore, if the function $p$ in the projective factor \(P=|y|p(r,s)\) is of the form \eqref{case3}, then the spray has constant curvature.

\medskip

\textbf{Case 4:} If $A\not\equiv0$ and $B\not\equiv0$, then we solve \eqref{ABeq2} and \eqref{ABeq3}.
Differentiating \eqref{ABeq2} with respect to \(u\) and multiplying
the result by \(A\), we get
\begin{equation*}
-2A^2BB'-AB^2A'+4A^2A'+2AA''=0.
\end{equation*}
From \eqref{ABeq3} we have $B'=\frac{1}{2}(B^3-3AB)$; substituting this into the above equation yields
\begin{equation*}
	-(AB^2)^2+(3A^2-A')AB^2+4A^2A'+2AA''=0.
\end{equation*}
From \eqref{ABeq2}, we obtain $AB^2=2A^2+2A'$; plugging this into
the above equation yields
\begin{equation*}
A^4-3(A')^2+AA''=0.
\end{equation*}
Moreover, since
\begin{equation*}
	(\frac{1}{A^2})''-2=-\frac{2(A^4-3(A')^2+AA'')}{A^4}=0,
\end{equation*}
it follows that
\begin{equation*}
\frac{1}{A^2}=u^2+2c_3u+c_4,
\end{equation*}
where $c_3,c_4$ are constants; hence
\begin{equation}\label{F1}
A=\frac{\delta}{\sqrt{u^2+2c_3u+c_4}},\quad A'=\frac{-\delta(u+c_3)}{(u^2+2c_3u+c_4)^{\frac{3}{2}}},
\end{equation}
where $\delta=\pm1$. From \eqref{ABeq2} we obtain $B=\eta\sqrt{2(A+\frac{A'}{A})}$, where $\eta=\pm1$. Plugging \eqref{F1} into the expression of $B$ yields
\begin{equation}\label{F2}
B=\eta\sqrt{\frac{2\left(\delta\sqrt{u^2+2c_3u+c_4}-u-c_3\right)}{u^2+2c_3u+c_4}}.
\end{equation}
On the domain where the square roots are well defined and
$B\not\equiv0$, a direct computation shows that \eqref{F1} and
\eqref{F2} satisfy \eqref{ABeq3} identically. Thus \eqref{ABeq2} and
\eqref{ABeq3} are both satisfied.
Plugging the expression of $A$ in \eqref{F1} and the expression of $B$ in \eqref{F2} back into \eqref{pre_u} yields
\begin{align}
\tilde{p}(u,v)&=
-\frac{\sqrt{2}\delta v+\eta\sqrt{\delta\sqrt{u^2+2c_3u+c_4}-u-c_3}}
{\sqrt{2}\left(\delta v^2+\sqrt{u^2+2c_3u+c_4}\right)+2\eta v\sqrt{\delta\sqrt{u^2+2c_3u+c_4}-u-c_3}} \notag \\
&=\frac{
\eta\left(\delta v^2-\sqrt{u^2+2c_3u+c_4}\right)\sqrt{\delta\sqrt{u^2+2c_3u+c_4}-u-c_3}-\sqrt{2}(u+v^2+c_3)v
}{
	\sqrt{2}\left(v^4+u^2+2c_3u+c_4+2v^2(u+c_3)\right)
}. \label{Case4p(u,v)_1}
\end{align}
Using \(p(r,s)=\tilde p(u,v)\), \((u,v)=(r-s^2,s)\), and
\eqref{Case4p(u,v)_1}, we obtain
\begin{align}
p(r,s)&=-\frac{\sqrt{2}\delta s+\eta\sqrt{\delta Q-r+s^2-c_3}}
{\sqrt{2}(\delta s^2+Q)+2\eta s\sqrt{\delta Q-r+s^2-c_3}} \label{sol4_1}\\
&=\frac{\eta(\delta s^{2}-Q)\sqrt{\delta Q-r+s^{2}-c_{3}}}{\sqrt2(r^2+2c_{3}r+c_{4})}
-\frac{(r+c_{3})s}{r^2+2c_{3}r+c_{4}}, \label{sol4_2}
\end{align}
where $\eta=\pm1$, $\delta=\pm1$, and $Q=\sqrt{(r-s^2)^2+2c_3(r-s^2)+c_4}$,  with $c_3$ and $c_4$ constants chosen such that $p(r,s)$ is well-defined, i.e.,
$\delta Q-r+s^2-c_{3}\ge 0,
(r-s^2)^2+2c_3(r-s^2)+c_4>0$, where the second inequality is strict---equivalently, $Q>0$---because $Q$ appears in the denominator of $B$ in \eqref{F2} and must be nonzero for $B$ to be defined.

Conversely, if $p(r,s)$ satisfies \eqref{sol4_2} (equivalently, \eqref{sol4}), then its first-order partial derivatives are
\begin{align*}
p_{r}&=\frac{s\left(r^2+2c_3r+2c_3^2-c_4\right)}
{\left(r^2+2c_3r+c_4\right)^2}
+
\frac{\sqrt{2}\eta\left(-3\delta Q-s^2+2r+2c_3\right)Q_r}
{4\left(r^2+2c_3r+c_4\right)\sqrt{\delta Q+s^2-r-c_3}}\\
&\quad+
\frac{
	\sqrt{2}\eta\left(\delta s^2-Q\right)
	\left(
	-4\delta(r+c_3)Q
	+3r^2+(-4s^2+6c_3)r-4c_3s^2+4c_3^2-c_4
	\right)}{4\left(r^2+2c_3r+c_4\right)^2\sqrt{\delta Q+s^2-r-c_3}
},
\\
p_s&=
-\frac{r+c_3}{r^2+2c_3r+c_4}-\frac{
	\eta\sqrt{2}\left[
	\left(3\delta Q+s^2-2r-2c_3\right)Q_s
	+
	\left(-2Q+\delta(-6s^2+4r+4c_3)\right)s
	\right]}{4\left(r^2+2c_3r+c_4\right)
	\sqrt{\delta Q+s^2-r-c_3}},
\end{align*}
where \[Q_r=\frac{r-s^2+c_3}{\sqrt{(r-s^2)^2+2c_3(r-s^2)+c_4}},\quad Q_s=\frac{2s(s^2-r-c_3)}{\sqrt{(r-s^2)^2+2c_3(r-s^2)+c_4}}.\]
Similarly, its second-order partial derivatives \(p_{ss}\), \(p_{rs}\), and \(p_{rr}\) can be calculated.
Substituting these derivatives into \eqref{const1} verifies that
\(p\) satisfies the system. Therefore, if the function \(p\) in the projective factor \(P=|y|p(r,s)\) of a projectively flat spherically symmetric spray is of the form given in \eqref{sol4_2}, then the spray has constant curvature.
\end{proofof}
\begin{remark}\label{remarkthm}
	The constants $c_3$ and $c_4$ in \eqref{sol4} can indeed be chosen so that the expression is well-defined. We give three such examples below; in each of them $\eta=\pm1$ is arbitrary.
	\begin{enumerate}[\rm (i)]
		\item Taking $\delta=1$, $c_3=0$, $c_4=1$, then $Q_1=\sqrt{(r-s^{2})^{2}+1}$,
		\[
		p(r,s)=\frac{\eta\,(s^{2}-Q_1)\sqrt{Q_1-r+s^{2}}}{\sqrt{2}\,(r^{2}+1)}-\frac{rs}{r^{2}+1}.
		\]
		Since $(r-s^{2})^{2}+1\geq1>0$, $Q_1-(r-s^{2})>0$ and $r^{2}+1\geq1>0$ for all $(r,s)$, this $p(r,s)$ is well-defined on the whole natural domain.
		\item Taking $\delta=1$, $c_3=-2$, $c_4=3$, then $Q_2=\sqrt{(r-s^{2})^{2}-4(r-s^{2})+3}$,
		\[
		p(r,s)=\frac{\eta\,(s^{2}-Q_2)\sqrt{Q_2-r+s^{2}+2}}{\sqrt{2}\,(r^{2}-4r+3)}-\frac{(r-2)s}{r^{2}-4r+3}.
		\]
		Writing $u=r-s^{2}$, one has $u^{2}-4u+3=(u-1)(u-3)>0$ for $0\leq u<1$, and $Q_2\geq u-2$ holds there trivially since $u-2<0$; hence this $p(r,s)$ is well-defined on the region $0\leq r-s^{2}<1$ with $r\neq1,3$.
		\item Taking $\delta=-1$, $c_3=-2$, $c_4=3$, then $Q_3=\sqrt{(r-s^{2})^{2}-4(r-s^{2})+3}$,
		\[
		p(r,s)=-\frac{\eta\,(s^{2}+Q_3)\sqrt{2-r+s^{2}-Q_3}}{\sqrt{2}\,(r^{2}-4r+3)}-\frac{(r-2)s}{r^{2}-4r+3}.
		\]
		With $u=r-s^{2}$ as above, $u^{2}-4u+3>0$ for $0\leq u<1$; moreover $2-u-Q_3\geq0$, since for $u\leq1$ both sides of $Q_3\leq2-u$ are nonnegative and squaring gives $u^{2}-4u+3\leq u^{2}-4u+4$. Hence this $p(r,s)$ is well-defined on the region $0\leq r-s^{2}<1$ with $r\neq1,3$.
	\end{enumerate}
\end{remark}

\section{Comparison with Finslerian Examples}\label{4}
In this section, we show that the sprays induced by the known projectively flat spherically symmetric Finsler metrics of constant flag curvature, namely the cases \(K=0\), \(K=-1\), and \(K=1\), are special instances of the solutions classified in Theorem~\ref{maintheorem1}. Recall that the projective factor is \(P(x,y)=|y|p(r,s)\), so it is determined by \(p(r,s)\). Each example follows the same pattern: we first recall the function \(p(r,s)\) associated with the metric, then identify the corresponding parameters in Theorem~\ref{maintheorem1}, and finally verify the equivalence by direct computation.

\begin{ex}
By Lemma \ref{classification of K}, a projectively flat spherically symmetric Finsler metric has constant flag curvature \(K=0\) if and only if $F$ is given by \eqref{K=0_F}. Computing the projective factor $P$ yields the corresponding function \(p_0=p(r,s)=P/\vert y\vert\),
\begin{equation}\label{p0_metric}
	p_0=0\quad or\quad p_0=\frac{s\pm\sqrt{c_0-r+s^2}}{c_0-r},
\end{equation}
where $c_0=1/c^2>0$ and $c$ is the constant in \eqref{K=0_F}. If $p_0=0$, this corresponds to \eqref{sol1}.
We now show that the second expression in \eqref{p0_metric} is a special case of \eqref{sol4}. To this end, we choose
\begin{equation}\label{c_K=0}
	c_3=-c_0,\qquad c_4=c_0^2,\qquad \delta=1,
\end{equation}
with $c_0>0$. On the domain $r-s^2\le c_0$, substituting these values into $Q$ in \eqref{sol4} yields
	\[
	Q=\sqrt{(r-s^2)^2+2c_3(r-s^2)+c_4}=\sqrt{(r-s^2)^2-2c_0(r-s^2)+c_0^2}=c_0-r+s^2,
	\]
and hence
	\[
	\delta Q-r+s^2-c_3 = (c_0-r+s^2)-r+s^2+c_0 = 2(c_0-r+s^2).
	\]
	Substituting these expressions into \eqref{sol4} gives
\begin{align*}
	p(r,s)&=
	\frac{\eta(\delta s^2-Q)\sqrt{\delta Q-r+s^2-c_3}}{\sqrt2(r^2+2c_3r+c_4)}
	-\frac{(r+c_3)s}{r^2+2c_3r+c_4}\\
	&=\frac{\eta(r-c_0)\sqrt{2(c_0-r+s^2)}-\sqrt2(r-c_0)s}{\sqrt2(r-c_0)^2}=\frac{\eta\sqrt{c_0-r+s^2}-s}{r-c_0}\\
	&=\frac{s-\eta\sqrt{c_0-r+s^2}}{c_0-r}.
\end{align*}

Since \(\eta=\pm1\), the resulting expression coincides with the second
expression in \eqref{p0_metric}. Therefore, with the parameters in
\eqref{c_K=0}, the solution \eqref{sol4} reproduces the metric-induced
function \(p_0(r,s)\); consequently, such sprays are induced by a
projectively flat spherically symmetric Finsler metric of constant
flag curvature \(K=0\).
\end{ex}
\begin{ex}
By Lemma \ref{classification of K}, a projectively flat spherically symmetric Finsler metric has constant flag curvature \(K=-1\) if and only if $F$ is given by \eqref{K=-1_F}. Computing the projective factor $P$ yields the corresponding function \(p_{-1}=p(r,s)=P/\vert y\vert\),
\begin{align*}
	p_{-1}&=\frac12\Bigg\{
	\frac{(c+1)^2s+(c+1)\sqrt{1-(c+1)^2(r-s^2)}}{1-(c+1)^2r}\nonumber\\
	&\quad+\frac{(c-1)^2s+(c-1)\sqrt{1-(c-1)^2(r-s^2)}}{1-(c-1)^2r}\Bigg\},
\end{align*}
where $c$ is a constant.
Here $p_{-1}$ is computed from the projective factor $P = F_{x^k} y^k / (2F)$, not read off directly from the expression of
$F$.
For subsequent computations, we use the
following equivalent form. Set
\[a=c+1,\qquad b=c-1,\qquad \gamma_a=\sqrt{1-a^2(r-s^2)},\qquad
\gamma_b=\sqrt{1-b^2(r-s^2)}.\]
Since \(1-a^2r = \gamma_a^2-a^2s^2=(\gamma_a-as)(\gamma_a+as)\), it follows that \[\frac{a^2s+a\gamma_a}{1-a^2r}
=\frac{a(as+\gamma_a)}{(\gamma_a-as)(\gamma_a+as)}
=\frac{a}{\gamma_a-as}.\]
The same computation applies to the \(\gamma_b\)-term. Hence
\begin{align}\label{p-1_metric}
p_{-1}=\frac12\left(\frac{a^2s+a\gamma_a}{1-a^2r}+\frac{b^2s+b\gamma_b}{1-b^2r}\right)=\frac12\left(\frac{a}{\gamma_a-as}+\frac{b}{\gamma_b-bs}\right).
\end{align}

We now verify that $p_{-1}$ matches one of the solutions in Theorem \ref{maintheorem1}, the case depending on $c$: if $c=0$, $p_{-1}$ reduces to \eqref{sol3} with $c_2=-1$; if $c=\pm1$, it reduces to \eqref{sol2} with $c_1=1/4$; if $c\neq0,\pm1$, it matches \eqref{sol4} with suitably chosen parameters.
\begin{enumerate}[\rm (i)]
\item
For $c=0$, we have $a=1,\ b=-1$, hence $\gamma_a=\gamma_b=\sqrt{1+s^2-r}$. From \eqref{p-1_metric},
\[
p_{-1}=\frac12\left(\frac1{\gamma_a-s}+\frac{-1}{\gamma_b+s}\right)=\frac{s}{1-r},
\]
which is precisely the special case $c_2=-1$ of \eqref{sol3} in Theorem \ref{maintheorem1}.
\item
For \(c=\pm1\), consider first \(c=1\), we have $a=2,\ b=0$. The term in \eqref{p-1_metric} involving $b$ vanishes identically, in this case $\gamma_a^2-4s^2=1-4r$, and hence
\[p_{-1}=\frac{1}{\gamma_a-2s}
=\frac{\gamma_a+2s}{\gamma_a^2-4s^2}=\frac{2s+\sqrt{1-4r+4s^2}}{1-4r}	=\frac{s+\sqrt{\frac14-r+s^2}}{2\times(\frac14-r)}.
\]
For $c=-1$, we have $a=0,b=-2$, a similar calculation yields
\[
p_{-1}=-\frac{1}{\gamma_b+2s}
=\frac{2s-\gamma_b}{\gamma_b^2-4s^2}=\frac{2s-\sqrt{1-4r+4s^2}}{1-4r}=\frac{s-\sqrt{\frac14-r+s^2}}{2\times(\frac14-r)}.
\]
Thus both cases coincide with the special case $c_1=1/4$ of \eqref{sol2} in Theorem \ref{maintheorem1}.
\item 		
For $c\neq0,\pm1$, we have $a,b\neq0$. We choose
\begin{equation}\label{c_K=-1}
c_3=-\frac{c^2+1}{(c^2-1)^2},\quad
c_4=\frac1{(c^2-1)^2},\quad
\delta=\operatorname{sgn}(c^2-1),\quad
\eta=-\operatorname{sgn}(c),
\end{equation}
in \eqref{sol4}; equivalently, in terms of \(a,b\),
\begin{equation}\label{tr_c_K=-1}
		c_3=-\frac{a^2+b^2}{2a^2b^2},\qquad
		c_4=\frac1{a^2b^2},\qquad
		\delta=\operatorname{sgn}(ab),\qquad
		\eta=-\operatorname{sgn}(a+b).
 		\end{equation}
Substituting these values into $Q$ in \eqref{sol4} yields
\begin{align*}
Q^2&=(r-s^2)^2+2c_3(r-s^2)+c_4=(r-s^2)^2-\frac{a^2+b^2}{a^2b^2}(r-s^2)+\frac1{a^2b^2}\\
&=\frac{\bigl(1-a^2(r-s^2)\bigr)\bigl(1-b^2(r-s^2)\bigr)}{a^2b^2}=\frac{\gamma_a^2\gamma_b^2}{a^2b^2},
\end{align*}
hence
		\begin{equation}\label{Q}
	 Q=\frac{\gamma_a\gamma_b}{\operatorname{sgn}(ab)ab}.
		\end{equation}
		Next, from the expressions for \(\gamma_a\) and \(\gamma_b\), we have $s^2-r=\frac12\left(\frac{\gamma_a^2-1}{a^2}+\frac{\gamma_b^2-1}{b^2}\right)$. Together with $\delta=\operatorname{sgn}(ab)$, this gives
		\[
		\begin{aligned}
\delta Q-r+s^2-c_3&=\operatorname{sgn}(ab)\times\frac{\gamma_a\gamma_b}{\operatorname{sgn}(ab)ab}+s^2-r+\frac{a^2+b^2}{2a^2b^2}\\
			&=\frac{b^2\gamma_a^2+a^2\gamma_b^2+2ab\gamma_a\gamma_b}{2a^2b^2}=\frac{(a\gamma_b+b\gamma_a)^2}{2a^2b^2},
		\end{aligned}	\]
and therefore
		\begin{equation}\label{sqrtQ}
		\sqrt{\delta Q-(r-s^2)-c_3}=\frac{\operatorname{sgn}(a+b)(a\gamma_b+b\gamma_a)}{\sqrt2\operatorname{sgn}(ab)ab}.
		\end{equation}
As shown in the proof of Theorem \ref{maintheorem1}, the solution \eqref{sol4} is equivalent to \eqref{sol4_1}, the two differing
only by rationalization. Since rationalization does not change the expression, we substitute the above results into \eqref{sol4_1} for computational convenience.
By \eqref{Q},\eqref{sqrtQ} and the chosen values of $\delta,\eta$ in \eqref{tr_c_K=-1}, we obtain
\[
\begin{aligned}
	p(r,s)&=
	-\frac{\sqrt{2}\delta s+\eta\sqrt{\delta Q-(r-s^2+c_3)}}
	{\sqrt{2}(\delta s^2+Q)+2\eta s\sqrt{\delta Q-(r-s^2+c_3)}}\\
	&=-\frac{
		\sqrt2\operatorname{sgn}(ab)s-\operatorname{sgn}(a+b)\times\dfrac{\operatorname{sgn}(a+b)(a\gamma_b+b\gamma_a)}{\sqrt2\operatorname{sgn}(ab)ab}}
	{\sqrt2\left(\operatorname{sgn}(ab)s^2+\dfrac{\gamma_a\gamma_b}{\operatorname{sgn}(ab)ab}\right)-2\operatorname{sgn}(a+b)s\times\dfrac{\operatorname{sgn}(a+b)(a\gamma_b+b\gamma_a)}{\sqrt2\operatorname{sgn}(ab)ab}}\\
	&=-\frac{
		2abs-(a\gamma_b+b\gamma_a)}
	{2\left(abs^2+\gamma_a\gamma_b\right)-2s(a\gamma_b+b\gamma_a)}=
	\frac{a\gamma_b+b\gamma_a-2abs}
	{2\bigl(abs^2-s(a\gamma_b+b\gamma_a)+\gamma_a\gamma_b\bigr)}\\
	&=\frac12\left(\frac{a}{\gamma_a-as}+\frac{b}{\gamma _b-bs}\right).
\end{aligned}
\]
This is exactly the expression \eqref{p-1_metric} for $p_{-1}$.
\end{enumerate}		

Therefore, for every value of $c$, the spray induced by such a projectively flat spherically symmetric Finsler metric coincides with one of the solutions in Theorem \ref{maintheorem1}; that is, the $K=-1$ metric-induced sprays are contained in our classification.
	\end{ex}

\begin{ex}
By Lemma \ref{classification of K}, a projectively flat spherically symmetric Finsler metric has constant flag curvature \(K=1\) if and only if $F$ is given by \eqref{K=1_F}. Since this form involves multi-valued
branches, the computation of geometric quantities becomes cumbersome. To obtain a concise real-valued expression, we use the equivalent real-parameterized construction in \cite{Shen}. Specifically,
for a projectively flat Finsler metric of constant flag curvature
\(K=1\) and for any angle \(\alpha\in[0,\pi/2)\), the
metric function \(F\) is expressed as
\begin{align}
F(x,y)&=\operatorname{Im}\left[\frac{-\langle x,y\rangle+i\sqrt{(e^{2i\alpha}+|x|^2)|y|^2-\langle x,y\rangle^2}}{e^{2i\alpha}+|x|^2}\right]\label{shen_K1}\\
&=\sqrt{\frac{\sqrt{\mathcal{A}}+\mathcal{B}}{2\mathcal{D}} + \left(\frac{\mathcal{C}}{\mathcal{D}}\right)^2} + \frac{\mathcal{C}}{\mathcal{D}}\label{K=1_F_alpha},
\end{align}
where
\begin{align*}
\mathcal{A}&= \left( \cos(2\alpha)|y|^2 + (|x|^2|y|^2 - \langle x,y \rangle^2) \right)^2 + \left( \sin(2\alpha)|y|^2 \right)^2,\\
\mathcal{B}&= \cos(2\alpha)|y|^2 + (|x|^2|y|^2 - \langle x,y \rangle^2),\\
\mathcal{C}&= \sin(2\alpha)\langle x,y \rangle,\qquad\mathcal{C}'=(\cos(2\alpha)+|x|^2)\langle x,y \rangle,\\
\mathcal{D}&= |x|^4 + 2\cos(2\alpha)|x|^2 + 1,
\end{align*}
and the projective factor $P$ is given by
\begin{equation}
P=\sqrt{\frac{\sqrt{\mathcal{A}}-\mathcal{B}}{2\mathcal{D}} - \left(\frac{\mathcal{C}}{\mathcal{D}}\right)^2} - \frac{\mathcal{C}'}{\mathcal{D}}.\label{shen_P}
\end{equation}

	We first note that setting $c=\tan\alpha$ in \eqref{K=1_F} gives \eqref{K=1_F_alpha}. To see the underlying equivalence, we compare the complex expressions inside \(\operatorname{Im}\) in \eqref{K=1_F} and \eqref{shen_K1}. Substituting the identities
	\[c+i=\frac{ie^{-i\alpha}}{\cos\alpha},\qquad(c+i)^2=-\frac{e^{-2i\alpha}}{\cos^2\alpha},\qquad(c+i)^4=\frac{e^{-4i\alpha}}{\cos^4\alpha},\]
	together with the rescaling
	\begin{equation}\label{rescale_k1}
		x=\frac{\tilde{x}}{\sqrt{c^2+1}},\qquad y=\frac{\tilde{y}}{\sqrt{c^2+1}}
	\end{equation}
	into \eqref{K=1_F}, we reduce \eqref{K=1_F} to \eqref{shen_K1} through algebraic simplification. Taking imaginary parts in each case gives the same real metric function. Consequently, \eqref{K=1_F} and \eqref{shen_K1} have the same structure and differ only by the rescaling of \(x\) and \(y\) in \eqref{rescale_k1}. In \cite{Shen}, the angle is taken in \(\alpha\in[0,\pi/2)\), which covers only the case \(c=\tan\alpha\ge 0\). For \(c<0\), the corresponding metric is obtained by choosing the other branch of the complex square root inside \(\operatorname{Im}\) in \eqref{shen_K1}. Hence the two formulas are equivalent for all \(c\in\mathbb{R}\).

Since \eqref{K=1_F} and \eqref{K=1_F_alpha} are equivalent, by \eqref{shen_P}, the corresponding function \(p_1=p(r,s)=P/\vert y\vert\) is given by
\begin{align}\label{p1_metric}
	p_1&=\frac{\sqrt{\left(\sqrt{(r-s^2+\cos(2\alpha))^2+\sin^2(2\alpha)}
			-r+s^2-\cos(2\alpha)\right)\left(r^2+2\cos(2\alpha)r+1\right)-2\sin^2(2\alpha)s^2}}{\sqrt2(	r^2+2\cos(2\alpha)r+1)}\nonumber\\
	&\quad-
	\frac{(r+\cos(2\alpha))s}
	{r^2+2\cos(2\alpha)r+1},
	\end{align}
	where \(\alpha\in[0,\pi/2)\). We choose
	\begin{equation}\label{c_K=1}
	c_3=\cos(2\alpha),\qquad c_4=1,\qquad \delta=1,\qquad \eta=-1,
	\end{equation}
	in \eqref{sol4}, and obtain
\begin{align*}
&Q=\sqrt{(r-s^2)^2+2c_3(r-s^2)+c_4}=\sqrt{(r-s^2+\cos(2\alpha))^2+\sin^2(2\alpha)}.
\end{align*}
On a small neighborhood of the origin in the \((r, s)\)-plane,
since
\(Q(0,0)=1\), we have \(Q-s^2>0\) on a sufficiently small
neighborhood of the origin. On this branch \(Q\ge s^2\), and with
the choices \(\delta=1\), \(\eta=-1\) in \eqref{c_K=1}, we have
\[
\eta(\delta s^2-Q)=Q-s^2=\sqrt{(Q-s^2)^2}.
\]
This yields
\begin{align}
	p(r,s)&=\frac{\eta(\delta s^{2}-Q)\sqrt{\delta Q-r+s^{2}-c_{3}}}{\sqrt2(r^2+2c_{3}r+c_{4})}
	-\frac{(r+c_{3})s}{r^2+2c_{3}r+c_{4}}\nonumber\\
	&=\frac{\sqrt{ \left(Q-r+s^{2}-\cos(2\alpha)\right)(Q-s^2)^2}}{\sqrt2(r^2+2\cos(2\alpha)r+1)}
	-\frac{(r+\cos(2\alpha))s}{r^2+2\cos(2\alpha)r+1}.\label{sim_prs}
\end{align}
On the other hand, we have
\begin{align*}	p_1=\frac{\sqrt{(Q-r+s^{2}-\cos(2\alpha))(r^2+2\cos(2\alpha)r+1)-2\sin^2(2\alpha)s^2}}{\sqrt2(r^2+2\cos(2\alpha)r+1)}-
\frac{(r+\cos(2\alpha))s}
{r^2+2\cos(2\alpha)r+1},
\end{align*}
and, since
\begin{align*}
\sin^2(2\alpha)=Q^2-(r-s^2+\cos(2\alpha))^2=\left(Q-r+s^2-\cos(2\alpha)\right)\times\left(Q+r-s^2+\cos(2\alpha)\right),
\end{align*}
it follows that
\begin{align}	p_1&=\frac{\sqrt{\left(Q-r+s^{2}-\cos(2\alpha)\right)\left(r^2+2\cos(2\alpha )r+1-2s^2(Q+r-s^{2}+\cos(2\alpha))\right)}}{\sqrt2(r^2+2\cos(2\alpha)r+1)}\nonumber\\
&\quad-\frac{(r+\cos(2\alpha))s}
{r^2+2\cos(2\alpha)r+1}\nonumber\\
&=\frac{\sqrt{\left(Q-r+s^{2}-\cos(2\alpha)\right)\left((r-s^2+\cos(2\alpha))^2+\sin^2(2\alpha)+s^4-2s^2Q\right)}}{\sqrt2(	r^2+2\cos(2\alpha)r+1)}\nonumber\\
&\quad-\frac{(r+\cos(2\alpha))s}
{r^2+2\cos(2\alpha)r+1}\nonumber\\
&=\frac{\sqrt{\left(Q-r+s^{2}-\cos(2\alpha)\right)(Q-s^2)^2}}{\sqrt2(	r^2+2\cos(2\alpha)r+1)}-
\frac{(r+\cos(2\alpha))s}
{r^2+2\cos(2\alpha)r+1}.\label{sim_p1}
\end{align}

Therefore, with the parameters in \eqref{c_K=1}, \eqref{sim_prs}
coincides with \eqref{sim_p1}, so the general expression \eqref{sol4}
reduces to \eqref{p1_metric}. This shows that the spray induced by a
projectively flat spherically symmetric Finsler metric with constant
flag curvature \(K=1\) is a special case of the classification
in Theorem~\ref{maintheorem1}.
\end{ex}

We conclude that the sprays induced by projectively flat metrics of constant flag curvature are all special cases in Theorem~\ref{maintheorem1}. Moreover, when the parameters take other values, one obtains sprays that are not induced by such metrics.

\section*{Declarations}
\noindent{\bf Conflict of interest.} There is no conflict of interest.

\noindent{\bf Availability of data and material.} Not applicable.

\noindent{\bf Declaration of Generative AI and AI-assisted technologies in the writing process.} AItools were used exclusively to assist with language polish and expression. The authors conducted all subsequent editing, verified the validity of the material, and assume sole responsibility for the final manuscript.

\end{document}